\documentclass[letter]{amsart}

\usepackage{amsmath}
\usepackage{amsthm}
\usepackage{amssymb}
\usepackage{url}
\usepackage[all,poly]{xy}
\SelectTips{cm}{}

\usepackage{diplomega_thm}
\usepackage{diplomega_math}

\makeatletter
\renewcommand{\@cite}[2]{%
  [#1\ifthenelse{\boolean{@tempswa}}{; #2}{}]}
\makeatother

\title[Measure Homology and Singular Homology]%
      {Measure Homology and Singular Homology are Isometrically Isomorphic}
      
\author{Clara L\"oh}
\address{Graduiertenkolleg \emph{Analytische Topologie und Metageometrie},
           Universit\"at M\"unster,  
           Ein\-stein\-str.~62,  
           48149 M\"unster, 
           Germany\\
           \url{http://www.math.uni-muenster.de/u/clara.loeh/}
        }
\email{clara.loeh@uni-muenster.de}%
\date{\today}

\begin{document}
\begin{abstract}
   Measure homology is a variation of singular homology designed by
   Thurston in his discussion of simplicial volume. Zastrow and Hansen
   showed independently that singular homology (with real
   coefficients) and measure homology coincide algebraically on the
   category of CW-complexes. It is the aim of this paper to prove
   that this isomorphism is \emph{isometric} with respect to the
   $\ell^1$-seminorm on singular homology and the seminorm on measure
   homology induced by the total variation. This, in particular,
   implies that one can calculate the simplicial volume via measure
   homology -- as already claimed by Thurston. For example, measure
   homology can be used to prove Gromov's proportionality principle of
   simplicial volume. 
\end{abstract}
\maketitle

\section{Introduction}

The simplicial volume is a topological invariant of oriented closed
connected manifolds, measuring the complexity of the fundamental
class. Despite its topological nature, the simplicial volume is linked
to Riemannian geometry in various ways~\cite{gromov}.

Originally, Gromov defined the simplicial volume to be the
$\ell^1$-seminorm of the fundamental class (with real
coefficients)\cite{munkholm}. In his famous lecture notes
\cite[Chapter~6]{thurston}, Thurston suggested an alternative
description of the simplicial volume for smooth manifolds: he replaced
singular homology and the $\ell^1$-seminorm by a new homology theory,
called smooth measure homology, and a corresponding
seminorm. However, except for the case of hyperbolic manifolds
\cite[Remark~0.1]{zastrow}, there is no published proof that these two
constructions result in the same simplicial volume. It is the purpose
of the present article to close this gap.

More generally, we prove the following:

\begin{satz}\label{isomisomthm}
  For all connected \cw-complexes~$X$, the inclusion $i_X
  \colon \csing * X \longrightarrow \mhch * X$ of the singular chain
  complex into the measure chain complex induces a natural isomorphism
  \[ \hsing * X \cong \hmh * X.\] This isomorphism is isometric with
  respect to the $\ell^1$-seminorm on singular homology and the
  seminorm on measure homology induced by the total variation.
\end{satz}

\begin{satz}\label{smisomisomthm}
  For all connected smooth manifolds~$M$ the canonical inclusions $\sminc_M
  \colon \smsing * M \longrightarrow \smmhch * M$ (of the smooth
  singular chain complex into the smooth measure chain complex) and
  $j_M \colon \smsing * M \longrightarrow \csing * M$ induce a natural
  isomorphism 
  \[ H_*(\sminc_M) \circ H_*(j_M)^{-1} 
     \colon \hsing * M \longrightarrow \hsmmh * M.\]
  This isomorphism is isometric with respect to the $\ell^1$-seminorm
  on singular homology and the seminorm on smooth measure homology
  induced by the total variation.
\end{satz}

We will now explain the occurring terminology in more detail.

The $\ell^1$-seminorm on singular homology with real coefficients is
the seminorm induced by the $\ell^1$-norm on the singular chain
complex: 

\begin{defi}
  Let $X$ be a topological space and let $k \in \N$. We define a
  norm~$\lone{\,\cdot\,}$ on the singular $k$-chains~$\csing k X$ with
  real coefficients by 
  \[ \lone c := \fsum_{\sigma \in \map[norm]{\ssim k} X} |a_\sigma| \]
  for all $c = \fsum_{\sigma \in \map[norm]{\ssim k}X} a_\sigma \cdot \sigma 
  \in \csing k X.$ This norm
  induces a seminorm on the $k$-th singular homology~$\hsing k X$ with
  real coefficients by
  \[ \lone \alpha := \inf \bigl\{ \lone c 
                          \bigm| c \in \csing k X, \bou(c) = 0, [c] = \alpha
                          \bigr\} \]
  for all $\alpha \in \hsing k X$.
\end{defi}

One of the most important properties of~$\lone{\,\cdot\,}$ is
``functoriality,'' i.e., for all continuous maps~$f$ the induced
homomorphism~$\hsing * f$ does not increase the seminorm. 
(Based on this property, Gromov introduced a more general framework of
functorial seminorms \cite[5.34]{gromovb}.)

Measure homology is a curious generalisation of singular homology: 
Let $X$ be a topological space, let $k \in \N$ and let $S_k(X) \subset
\map[norm]{\ssim k}{X}$ be some set of singular simplices. The idea of
measure homology is to think of a singular chain $\fsum_{\sigma \in
S_k(X)} a_\sigma \cdot \sigma$ with real coefficients as a signed measure
on~$S_k(X)$ having the mass~$a_\sigma$ on the set~$\{\sigma \}.$ The
measure chain complex consists of all signed measures on~$S_k(X)$
satisfying some finiteness condition. 

Thus the measure chain complex is larger than the singular chain
complex and hence provides more room for constructions such as
``smearing'' \cite[page~6.8ff, page~547ff]{thurston,ratcliffe}. The
other side of the coin is that it is quite hard to gain a geometric
intuition of more complicated measure chains.

Depending on the choice of the mapping spaces~$S_k(X)$ and their
topology, there are two main flavours of measure homology:

\begin{itemize}
  \item One for general topological spaces using the
        compact-open topology on the set of all singular
        simplices -- the corresponding chain complexes and homology
        groups are denoted by~$\mhch * X$ and~$\hmh * X$.
  \item One for smooth manifolds using the \cone-topology
        on the set of smooth singular simplices -- the corresponding
        chain complexes and homology groups are denoted by~$\smmhch *
        X$ and~$\hsmmh * X$.
\end{itemize}

In both cases, measure homology is equipped with the seminorm induced
by the total variation on the chain level. 

Measure homology of the second kind (so-called smooth measure
homology) was introduced by Thurston~\cite[page~6.6]{thurston}. Some
basic properties of smooth measure homology are also listed in
Ratcliffe's book~\cite[\S11.5]{ratcliffe}. A thorough treatment of
measure homology for general spaces is given in the papers of Zastrow
and Hansen \cite{zastrow,hansen}. In both papers it is shown that
measure homology for \cw-complexes coincides \emph{algebraically} with
singular homology with real coefficients.

The general idea of the proof of Theorem~\ref{isomisomthm}
and~\ref{smisomisomthm} is to take a dual point of view: singular
homology and the $\ell^1$-seminorm admit a dual concept, called
bounded cohomology. The key property is that the canonical seminorm on
bounded cohomology and the $\ell^1$-seminorm are intertwined by a
duality principle (Theorem~\ref{duality}).

Unlike the $\ell^1$-seminorm, bounded cohomology and its seminorm
are rather well understood by Ivanov's homological algebraic approach
\cite{ivanov,monod}. In our proof of Theorem~\ref{isomisomthm}
and~\ref{smisomisomthm}, we take advantage of a special cochain
complex~$\bfch * X$ computing bounded cohomology. More precisely, we
construct a dual~$\bmcch * X$ of the measure chain complex (together
with a corresponding duality principle) fitting into a commutative 
diagram
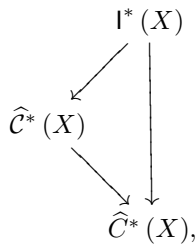
\begin{align*}
  \xymatrix@=3em@dr{
    \bfch  * X \ar[dr] \ar[d] & \\
    \bmcch * X \ar[r]         & \bchn * X,
  }
\end{align*}
where $\bchn * X$ is the cochain complex defining bounded
cohomology. The crux of this diagram is that the vertical arrow induces an
isometric isomorphism on the level of cohomology and that the arrows
on the left do not increase the seminorm. Then the duality principles
allow us to deduce that the algebraic isomorphism~$\hsing * X \cong
\hmh * X$ must be isometric. 

For simplicity, we only prove Theorem~\ref{isomisomthm} in detail.
The smooth case, requiring a small detour to smooth singular homology,
is considered briefly in Section~\ref{smoothsec} (and is also
explained in the author's diploma thesis \cite{strohm}).

This paper is organised as follows: Measure homology (the non-smooth
version) is defined in Section~\ref{mhsec}. Section~\ref{dualsec} is
concerned with the dual point of view, i.e., the construction of a
dual for measure homology and the derivation of a corresponding
duality principle.  A proof of Theorem~\ref{isomisomthm} is presented
in Section~\ref{proofsec}. In Section~\ref{smoothsec}, we have a
glimpse at the smooth universe, that is at smooth measure homology and
at a proof of Theorem~\ref{smisomisomthm}. Finally, in
Section~\ref{applications}, we list some applications to the
simplicial volume, including Gromov's proportionality principle and
some of its consequences.

\section{Measure homology}\label{mhsec}

In this section, our basic object of study, measure homology, is
introduced. In Section~\ref{mhalgiso}, we describe the algebraic
isomorphism between singular homology and measure homology. The smooth
case is deferred to Section~\ref{smoothsec}.

\subsection{Definition of measure homology}

Before stating the precise definition of measure homology, we recall
some basics from measure theory:

\begin{defi}
  \newcommand{\sigb}{{\mathcal B}}%
  Let $(X, \siga)$ be a measurable space. 
  \begin{itemize}
    \item A map $\mu \colon \siga \longrightarrow \R \cup \{ \infty, -\infty\}$
          is called a \defin{signed measure} if $\mu(\emptyset)=0$,
          not both $\infty$ and $-\infty$ are contained in the image
          of~$\mu$, and $\mu$ is $\sigma$-additive. 
    \item A \defin{null set} of a signed measure~$\mu$ on~$(X,\siga)$ is a
          measurable set~$A \in \siga$ with $\mu(B) = 0$ for all $B\in
          \siga$ with $B \subset A$.
    \item A \defin{determination set} of a signed measure~$\mu$
          on~$(X,\siga)$ is a subset~$D$ of~$X$ such that each
          measurable set contained in the complement of~$D$ is a
          $\mu$-null set.
    \item The \defin{total variation} of a signed measure~$\mu$
          on~$(X,\siga)$ is given by
          \[ \tov \mu := \sup_{A \in \siga} \mu(A) 
                       - \inf_{A \in \siga} \mu(A). 
          \]    
    \item For $x \in X$ the \defin{atomic measure} concentrated in~$x$
          is denoted by~$\delta_x$. 
    \item If $f \colon (X,\siga) \longrightarrow (Y,\sigb)$ is a
          measurable map and $\mu$ is a signed measure on~$(X,\siga),$
          then 
          \[ \fa{B \in \sigb} \mu^f(B) := \mu\bigl(f^{-1}(B)\bigr)\]
          defines a signed measure~$\mu^f$ on~$(Y,\sigb)$. \qedhere  
  \end{itemize}
\end{defi}

As indicated in the introduction, the measure homology chain complex
consists of measures that respect some finiteness condition on the set
of all singular simplices. 

\begin{defi}
  Let $X$ be a topological space and let $k \in \N$. 
  \begin{itemize}
    \item The $k$\defin{-th measure chain group}, denoted by~$\mhch k
          X$, is the $\R$-vector space of signed measures 
          on~$\map[norm]{\ssim k}X$ possessing a
          compact determination set and finite total variation.  
          Here $\map[norm]{\ssim k} X$ is equipped with the
          compact-open topology and the corresponding Borel
          $\sigma$-algebra. The elements of~$\mhch k X$ are called 
          \defin{measure} $k$\defin{-chains}.   
    \item For each $j \in \{ 0, \dots, k+1\}$ the inclusion~$\bou_j
          \colon \ssim k \longrightarrow \ssim{k+1}$ of the $j$-th
          face induces a continuous map $\map[norm]{\ssim {k+1}} X
          \longrightarrow \map[norm]{\ssim k}X$ and hence a homomorphism 
          (which we will also denote by~$\bou_j$)
          \begin{align*}
             \bou_j \colon \mhch{k+1} X & \longrightarrow \mhch{k} X\\
                                   \mu & \longmapsto
                    \mu^{(\sigma \mapsto \sigma \circ \bou_j)}.
          \end{align*} 
          The \defin{boundary operator} of measure chains is then
          defined by
          \begin{align*}
             \bou := \sum_{j=0}^{k+1} (-1)^j \cdot \bou_j
             \colon \mhch{k+1} X & \longrightarrow \mhch k X.
          \end{align*}
    \item The $\R$-vector space
          $\hmh k X := H_k\bigl( \mhch * X, \bou \bigr)$
          is called the $k$\defin{-th measure homology group} of~$X$.
    \item The total variation $\tov{\,\cdot\,}$ turns~$\mhch k X$ into
          a normed vector space and thus   
          induces a seminorm on~$\hmh k X$ as follows: For 
          all~$\mu \in \hmh k X$ we define
          \[ \mhnorm \mu := \inf \bigl\{ \tov \nu
                                 \bigm|  \nu \in \mhch k X,
                                         \bou(\nu) = 0,
                                         [\nu] = \mu
                                 \bigr\}. \qedhere\]
  \end{itemize}
\end{defi}

Zastrow showed that $(\mhch * X, \bou)$ indeed is a chain complex
\cite[Corollary~2.9]{zastrow}. Hence measure homology is well-defined. 
Each continuous map $f \colon X \longrightarrow Y$ induces a chain
map \cite[Lemma-Definition~2.10(iv)]{zastrow}
  \begin{align*}
     \mhch * f \colon \mhch * X & \longrightarrow \mhch * Y \\
                            \mu & \longmapsto     \mu^f
  \end{align*}
which obviously does not increase the total variation. Therefore, we obtain a
homomorphism $\hmh * f \colon \hmh * X \longrightarrow \hmh * Y$
satisfying $\mhnorm{\hmh * f (\mu)} \leq \mhnorm{\mu}$ for all~$\mu
\in \hmh * X$. 
Clearly, this turns ${\mathcal H}_*$ into a functor. Moreover, the
functor~${\mathcal H}_*$ is homotopy invariant
\cite[Lemma-Definition~2.10(vi)]{zastrow}. 

Analogously to singular homology, relative measure homology groups can
be defined \cite[Lemma-Definition~2.10(ii)]{zastrow}.

\subsection{The algebraic isomorphism}\label{mhalgiso}

There is an obvious norm-preserving inclusion of the singular chain
complex into the measure chain complex: 

\begin{defi}
  If $X$ is a topological space and $k \in \N$, we write 
  \begin{align*}
    i_X \colon \csing k X & \longrightarrow \mhch k X\\ 
    \fsum_{\sigma \in \map[norm]{\ssim k}{X}} 
                            a_\sigma \cdot \sigma
                          & \longmapsto
    \fsum_{\sigma \in \map[norm]{\ssim k}{X}} 
                            a_\sigma \cdot\delta_\sigma. 
    \qedhere
  \end{align*}
\end{defi}

\begin{qedbem}
  This inclusion induces a natural chain map $\csing * X
  \longrightarrow \mhch * X$ (also denoted by $i_X$) which is
  norm-preserving.  
\end{qedbem}

Establishing the Eilenberg-Steenrod axioms for measure homology,
Hansen~\cite{hansen} and Zastrow~\cite{zastrow} independently proved
the following theorem:

\begin{satz}\label{algebraicisomthm}
  For all \cw-complexes~$X$, the inclusion $i_X \colon \csing * X
  \longrightarrow \mhch * X$ induces a natural isomorphism (of real
  vector spaces)
  \[ \hsing * X \cong \hmh * X.\]
\end{satz}

However, there are spaces for which singular homology and measure
homology do not coincide \cite[Section~6]{zastrow}.

In Section~\ref{proofsec}, we prove the main result, namely that
these isomorphisms are compatible with the induced seminorms on~$H_*$
and~${\mathcal H}_*$. 
Consequences of this theorem and of its smooth analogue
(Theorem~\ref{smisomisomthm}) are discussed in Section~\ref{applications}.

\section{A dual point of view}\label{dualsec}

Rather than attempting to investigate the functorial
seminorms~$\lone{\,\cdot\,}$ and~$\mhnorm{\,\cdot\,}$ on singular
homology and measure homology directly on the chain level, we take a
dual point of view: we make use of bounded cohomology and the duality
principle (Theorem~\ref{duality}) to
compute~$\lone{\,\cdot\,}$. Analogously, we also construct a dual for
measure homology and derive a corresponding duality principle
(Section~\ref{mhdualitysec}).  In Section~\ref{specialsec}, a special
cochain complex for bounded cohomology is introduced, which turns out
to be very convenient in our setting.

\subsection{Bounded cohomology}

Bounded cohomology is the functional analytic twin of singular
cohomology. It is constructed via the topological dual of the singular
chain complex instead of the algebraic one.  The corresponding norm
for singular cochains is therefore the supremum norm:

\begin{defi}
  Let $X$ be a topological space and $k \in \N.$ For a singular
  cochain $f \in \cosing{k}{X}$ the (possibly infinite)  
  \defin{supremum norm} is defined by 
  \[ \supn{f} := \sup\bigl\{ |f(\sigma)|
                      \bigm| \sigma \in \map[norm]{\ssim k}{X}
                      \bigr\} .\]
  This induces a seminorm on~$\cohsing{k}{X}$ by
  \[ \supn{\varphi} := 
     \inf \bigl\{ \supn{f} \;\big|\; 
             f \in \cosing{k}{X}, \cobou(f) =0, [f] = \varphi
          \bigr\}
   \]
  for all $\varphi \in \cohsing{k}{X}$. We write
  $ \bchn{k}{X} 
     := \bigl\{ f \in \cosing{k}{X} \;\big|\; \supn{f} < \infty 
        \bigr\}
  $ 
  for the vector space of \defin{bounded }$k$\defin{-cochains}.
\end{defi}

It is easy to see that the coboundary operator~$\cobou$ on the singular cochain
complex~$\cosing{k}{X}$ satisfies   
$\cobou\bigl( \bchn{k}{X}\bigr) \subset \bchn{k+1}{X}.$
Thus $\bchn{*}{X}$ is a cochain complex.

\begin{defi}\label{bcdef}
  Let $X$ be a topological space and $k \in \N.$ 
  \begin{itemize}
    \item The $k$\defin{-th bounded cohomology group of}~$X$ is defined by
          \[\bch k X := H^k\bigl( \bchn * X,\cobou|_{\bchn * X}\bigr).\]
    \item The supremum norm on~$\bchn{k}{X}$ induces a seminorm on~$\bch{k}{X}$
          by
          \[ \supn{\varphi} := \inf \bigl\{  \supn{f} \;
                                     \big|\; f  \in \bchn{k}{X}, \cobou(f) = 0,
                                            [f] = \varphi
                                     \bigr\}.\]
          for all~$\varphi \in \bch k X$. \qedhere 
  \end{itemize}
\end{defi}

Overviews of bounded cohomology (of spaces) are given by Ivanov
\cite{ivanov}, Gromov~\cite{gromov}, and Brooks \cite{brooks}, 
where also the more peculiar aspects of bounded cohomology are
explained.  For example, bounded cohomology depends only on the
fundamental group \cite[Corollary~(A) on page~40,
Theorem~(4.1)]{gromov,ivanov} and does not satisfy the excision axiom
\cite[\S3(a),\S5]{brooks,mitsumatsu}.

\subsection{Duality}

The duality principle (Theorem~\ref{duality}) shows an important
aspect of bounded cohomology: bounded cohomology can be used to
compute the  $\ell^1$-seminorm (and hence the simplicial
volume). Since bounded cohomology is much better
understood (in view of the techniques presented in
Ivanov's paper~\cite{ivanov}) than the seminorm on homology, duality
leads to interesting applications. For example, bounded cohomology can
be used to give estimates for the simplicial volume of products and
connected sums of manifolds \cite[page~10]{gromov}. 

Moreover, duality plays a central r\^ole in the proof that measure
homology and singular homology are \emph{isometrically} isomorphic
(see Section~\ref{proofsec}). 

\begin{satz}[Duality Principle]\label{duality}
  Let $X$ be a topological space, $k\in \N$ and $\alpha \in
  \hsing{k}{X}.$
  \begin{enumerate}
    \item Then $\lone{\alpha} = 0$ if and only if
          \[ \fa{\varphi \in \bch{k}{X}} \krp\varphi\alpha =0.\]
    \item If $\lone{\alpha} > 0,$ then
          \[ \lone\alpha 
             = \sup\Bigl\{ \frac{1}{\supn{\varphi}} 
                   \Bigm| \text{$\varphi\in\bch{k}{X}$,
                                $\krp\varphi\alpha = 1$}
                   \Bigr\}.\]
  \end{enumerate}
\end{satz}

The angle brackets~$\krp{\,\cdot\,}{\cdot\,}$ refer to the Kronecker
product on~$\bch * X \otimes \hsing * X$ defined by
evaluation, just as the ordinary Kronecker product on~$\cohsing * X
\otimes \hsing * X$.

This duality was discovered by Gromov \cite{gromov}.
A detailed proof -- based on the Hahn-Banach theorem
-- is given in the book of Benedetti and Petronio
\cite[Proposition~F.2.2]{bp}.

\subsection{A special cochain complex}\label{specialsec}

Bounded cohomology, as defined in Definition~\ref{bcdef}, is quite
hard to calculate. However, Ivanov found a homological algebraic route
to bounded cohomology via strong relatively injective resolutions of
Banach modules \cite{ivanov,monod}. In Section~\ref{proofsec}, the
resolution discussed in the following paragraphs saves the day:

\begin{defi}\label{specialdef}
  Let $X$ be an arcwise connected space with universal covering~$\ucov
  X$, and let $k \in \N$.
  \begin{itemize}
    \item Then $\pi_1 (X)$ acts from the left on the vector 
          space~$\map[norm]{\ucov X^{k+1}}{\R}$ of continuous functions
          $\ucov X ^{k+1} \longrightarrow \R$ by
          \[ g \cdot f := \bigl( (x_0, \dots, x_k) \mapsto
                                 f(x_0 \cdot g,  \dots, x_k \cdot g)
                          \bigr) \]
          for all $f \in \map[norm]{\ucov X ^{k+1}}\R$ and all $g \in
          \pi_1(X)$.
    \item The subset of bounded functions
          in~$\map[norm]{\ucov X^{k+1}}\R$ is denoted by~$\bbch k X$ 
          and we use the abbreviation
          $\bfch k X$ for the functions in~$\bbch k X$ that are invariant  
          under the above $\pi_1(X)$-action. 
    \qedhere
  \end{itemize}
\end{defi}

How can we turn $\bbch * X$ into a cochain complex? The vector space
$\bchn[norm] k {\ucov X}$ can be identified with the space of bounded
functions $\map[norm]{\ssim k}{\ucov X} \longrightarrow \R$ (under this
identification, the norm $\supn{\, \cdot\,}$ just becomes the supremum
norm). 
Now $\bbch k X$ can be viewed as a subspace of~$\bchn[norm] k {\ucov X}$,
namely as the space of those bounded functions $\map[norm]{\ssim k}{\ucov X}
\longrightarrow \R$ that only depend on the vertices of the simplices
(and are continuous in the vertices). In Gromov's terminology those
functions would be called ``straight bounded continuous cochains''
\cite[Section~2.3]{gromov}, inspiring the sans-serif notation. 

It is clear that the coboundary operator on~$\bchn[norm] * {\ucov X}$
restricts to $\bbch * X$ and that the operations of $\pi_1(X)$ on
$\bbch * X$ and
$\bchn[norm] * {\ucov X}$ are compatible with the above inclusion
map. This makes $\bbch * X$ a subcomplex of~$\bchn[norm] * {\ucov
X}$. In other words, the homomorphism
\begin{align*}
  u^k \colon \bbch k X & \longrightarrow \bchn[norm] k {\ucov X} \\
                     f & \longmapsto
                         \bigl( \sigma \mapsto
                                f(\sigma(e_0), \dots, \sigma(e_k))
                         \bigr)
\end{align*}
yields a $\pi_1(X)$-equivariant cochain map. Here, $e_0, \dots, e_k$
are the vertices of~$\ssim k$. 

\begin{satz}\label{specialthm}%
  Let $X$ be a connected locally finite \cw-complex and let $p \colon \ucov
  X \longrightarrow X$ be its universal covering. 
  For
  $\sigma \in \map[norm]{\ssim k} X$ we denote by~$\ucov \sigma \in
  \map[norm]{\ssim k}{\ucov X}$ some lift of~$\sigma$ with respect
  to~$p$. 
  The cochain map~$v \colon \bfch * X \longrightarrow \bchn * X$ given by
  \begin{align*} 
    \bfch k X & \longrightarrow \bchn k X \\
            f & \longmapsto     \bigl( \sigma \mapsto
                                       f(\ucov\sigma(e_0),\dots,
                                         \ucov\sigma(e_k)
                                        )
                                \bigr)
  \end{align*}
  induces an isometric isomorphism $H^*\bigl(\bfch * X\bigr) \cong \bch * X.$  
\end{satz}

Covering theory shows that $v$ is a well-defined cochain map.

\begin{qedbem}
  Since all connected smooth manifolds are triangulable (and hence
  locally finite \cw-complexes), the theorem applies in particular to
  this case.
\end{qedbem}

\begin{bew}
  The universal covering map~$p$ induces an isometric isomorphism 
  \[ \widehat p ^* \colon \bch * X \longrightarrow H^*\bigl( \bchn[norm] *
  {\ucov X}^{\pi_1(X)}\bigr)\] 
  where $\bchn[norm] * {\ucov X}^{\pi_1(X)}$
  denotes the subcomplex of $\pi_1(X)$-fixed points under the canonical $\pi_1(X)$-action
  on~$\map[norm]{\ssim k}{\ucov X}$ induced by the $\pi_1(X)$-action
  on~$\ucov X$ \cite[proof of Theorem~(4.1)]{ivanov}. 
  By construction, the triangle 
  \begin{align*}
    \xymatrix@R=2.5em@C=0.7em{
      H^*\bigl( \bfch * X\bigr) \ar[rr]^{H^*(v)} \ar[rd]_{H^*(u|)}
      & & \bch * X \ar[dl]^{\widehat p^*} \\
       & H^*\bigl( \bchn[norm]*{\ucov X}^{\pi_1(X)} \bigr)& 
    }
  \end{align*}
  is commutative (where $u|$ denotes the restriction of~$u$ to the
  $\pi_1(X)$-fixed points). 
  The fact that there exists a canonical 
  isometric isomorphism 
  \[ H^*\bigl(\bfch * X\bigr)
     \cong H^*\bigl( \bchn[norm] * {\ucov X}^{\pi_1(X)}\bigr)
  \] 
  follows from work of 
  Monod \cite[Theorem~7.4.5]{monod}. The homomorphism~ $u$ fits into
  the ladder
  \begin{align*}
    \xymatrix@=2em{
        0  \ar[r]             
      & \R \ar[r]^-\varepsilon \ar[d]_{\id}
      & \bbch 0X \ar[r]        \ar[d]_{u^0}
      & \bbch 1X \ar[r]        \ar[d]_{u^1}
      & \dots \\
     0 \ar[r] & \R \ar[r]_-\varepsilon
              & \bchn[norm]0{\ucov X} \ar[r] 
              & \bchn[norm]1{\ucov X} \ar[r]
              & \cdots\text{\makebox[0pt]{\phantom{(}}}\smash{,}
    }
  \end{align*}
  whose rows are strong resolutions of the trivial Banach
  $\pi_1(X)$-module~$\R$ by relatively injective $\pi_1(X)$-modules
  \cite[Theorem~7.4.5, proof of Theorem~(4.1)]{monod,ivanov}. 
  Hence we obtain that the induced map $H^*(u|)$ must be this 
  isometric isomorphism \cite[Lemma~7.2.6]{monod}. Therefore, 
  the composition $H^*(v) = (\widehat p^*)^{-1} \circ H^*(u|)$ is also
  an isometric isomorphism. 
\end{bew}

\subsection{A dual for measure homology}\label{mhdualitysec}

In order to develop a duality principle in the setting of measure
homology, we have to construct a ``dual''~$\smash{\bmc * X}$ playing the
r\^ole of the bounded cohomology groups~$\smash{\bch * X}$ in the singular
theory:

If $c = \fsum_{\sigma \in \map[norm]{\ssim k}X} a_\sigma \cdot \sigma \in
\csing k X$ 
is a singular chain and $f \in \bchn k X$ is a singular cochain,
their Kronecker product is given by
\[    \krp{f}{c} = f(c)
    = \fsum_{\sigma \in \map[norm]{\ssim k} X} a_\sigma \cdot
      f(\sigma)
    \in \R. \]
If we think of $c$ as a linear combination of atomic measures, this
looks like an integration of~$f$ over~$c$. Hence our ``dual'' in
measure homology consists of (bounded) functions that can be
integrated over measure chains:

\begin{defi}
  Let $X$ be a topological space and let $k\in\N$. 
  We define 
  \[  \bmcch k X
      :=  \bigl\{ f \colon \map[norm]{\ssim k}{X} \longrightarrow \R
          \bigm | \text{$f$ is Borel measurable and bounded}
          \bigr\} \]
  and (where $f(\bou(\sigma))$ is an abbreviation for~$\sum_{j =0}^{k+1} 
  (-1)^j \cdot f(\face j \sigma)$)
  \begin{align*}
     \cobou \colon \bmcch k X & \longrightarrow  \bmcch {k+1} X\\
                            f & \longmapsto      
                            \bigl( \sigma \mapsto (-1)^{k+1} \cdot f(\bou(\sigma
)) 
                            \bigr). \qedhere
  \end{align*}
\end{defi}

It is not hard to see that this map $\cobou \colon \bmcch k X
\longrightarrow \bmcch {k+1} X$ is indeed well-defined and that it
turns $\bmcch * X$ into a cochain complex.

\begin{defi}
  Let $X$ be a topological space and let $k \in \N$. The 
  $k$\defin{-th bounded measure cohomology group of}~$X$ is given by
  $ \bmc k X := H^k \bigl(\bmcch * X, \cobou \bigr).$  
  We write 
  $\supn{\,\cdot\,}$ for the seminorm 
  on~$\bmc k X$ which is induced
  by the supremum norm on~$\bmcch k X$.
\end{defi}

As a second step, we have to generalise the Kronecker product to
bounded measure cohomology. As indicated above, our Kronecker product
relies on integration:

\begin{defi}
  Let $X$ be a topological space and let $k \in
  \N$. The \defin{Kronecker product} 
  of $\mu \in \mhch k X$ and~$f \in \bmcch k X$ is defined as
  \[ \krp f \mu := \int f \;d\mu. \]
  If $\mu$ is a measure cycle and $f$ is a cocycle, we write
  $ \krp[big]{[f]}{[\mu]} := \krp f \mu = \int f \;d\mu. $
\end{defi}

The integral is defined and finite, since the elements of $\mhch k X$
are (signed) measures of finite total variation and the elements
of~$\bmcch k X$ are bounded measurable functions.  Moreover, the
integral is obviously bilinear. The transformation formula implies
that the definition of the Kronecker product on (co)homology does not
depend on the chosen representatives.  Hence the Kronecker product is
well-defined and bilinear.

\begin{lem}\label{mhkrp} 
  Let $X$ be a topological space and let~$k \in \N$. 
  The Kronecker product defined above is compatible with the
  Kronecker product on bounded cohomology 
  in the following sense: 
  for all $ f \in \bmcch k X$ and all $ c \in \csing k X$,
  \[ \krp {v_2(f)} c = \krp f {i_X (c)}, \]
  where $v_2(f)$ denotes the linear extension of~$f \colon
  \map[norm]{\ssim k} X \longrightarrow \R$ to the vector 
  space~$\csing k X.$          
  Passage to (co)homology yields for all $\varphi \in \bmc k
  X$ and all $\alpha \in \hsing k X$ 
  \[ \krp[big]{H^k(v_2)(\varphi)}{\alpha}
     = \krp[big]{\varphi}{H_k(i_X)(\alpha)}. 
  \]
\end{lem}
\begin{bew}
  Since both Kronecker products are bilinear, it suffices to
  consider the case where $c$ consists of a single singular 
  simplex~$\sigma \in \map[norm]{\ssim k} X$. Then the left hand side -- by
  definition -- evaluates to~$f(\sigma)$. Since $i_X(c)=
  i_X(\sigma) = \delta_\sigma$ is the atomic measure on~$\map[norm]
  {\ssim k} X$
  concentrated in~$\sigma$, we obtain 
  \[ \krp[big] f {i_X(c)} = \int f \;d\delta_\sigma = 1 \cdot f(\sigma)
  \]
  for the right hand side. The corresponding equality in (co)homology
  follows because $v_2 \colon \bmcch * X \longrightarrow \bchn *
  X$ is easily recognised to be a chain map.
\end{bew}

The above Kronecker product leads to the following (slightly weakened)
duality principle: 

\begin{lem}[Duality Principle of Measure Homology]\label{mhduality}%
  Let $X$ be a topological space, let $k\in\N$, and $\alpha
  \in \hmh k X$.
  \begin{enumerate}
    \item If $\mhnorm \alpha = 0$, then $\krp \varphi \alpha = 0$ 
          for all $\varphi \in \bmc k X$.
    \item If $\mhnorm \alpha > 0$, then 
          \[ \mhnorm \alpha \geq \sup 
                  \Bigl\{ \frac 1 {\supn \varphi} 
                  \Bigm|  \text{$\varphi \in \bmc k X$, $\krp
                          \varphi \alpha = 1$}
                  \Bigr\}.\]
  \end{enumerate}
\end{lem}
\begin{bew}
  Let $\varphi \in \bmcch k X$. Assume that $\mu \in \mhch k
  X$ is a measure cycle representing~$\alpha$ and $f \in \bmcch k X$
  is a cocycle representing~$\varphi$. 
  If $\krp \varphi \alpha =1$, then 
  \begin{align*}
     1   =    \lvert\krp \varphi \alpha\rvert   
         =    \Bigl|\int f \; d\mu\Bigr|
       & \leq  \supn f \cdot \tov\mu. 
  \end{align*}
  Taking the infimum over all representatives results in
  $1 \leq \supn \varphi \cdot \mhnorm \alpha.$
  In particular, if there exists such a $\varphi$, then
  \[ \mhnorm \alpha \geq \frac 1 {\supn \varphi} > 0.\]
  Now the lemma is an easy consequence of this inequality.
\end{bew}

\begin{bem}
  A posteriori we will be able to conclude -- in view of
  Theorem~\ref{isomisomthm}, Theorem~\ref{duality} and Lemma~\ref{mhkrp} -- 
  that in the first part of the lemma ``if and only if'' is also true
  and that in the second part equality holds.
\end{bem}

\section{Proving the isometry}\label{proofsec}

This section is devoted to the proof of Theorem~\ref{isomisomthm}. 
To show that the algebraic isomorphism between singular homology and
measure homology is isometric, we proceed in two steps. First, we
prove the theorem in the special case of connected locally finite
\cw-complexes. In the second step we generalise this result using a
colimit argument.

\subsection{First step -- connected locally finite \cw-complexes}\label{firststep}

We investigate the dual~$\bmcch * X$ by means of the
complex~$\bfch * X$ introduced in Definition~\ref{specialdef}: recall
that the vector space $\bfch k X$ is the set of all bounded functions
in~$\map[norm]{\ucov X^{k+1}}{\R}$ that are $\pi_1(X)$-invariant.
Then the key to the proof of Theorem~\ref{isomisomthm} is a careful
analysis of the diagram 
\begin{align*}
  \xymatrix@=3em@dr{
    \bfch  * X \ar^{v}[dr] \ar_{v_1}[d] & \\
    \bmcch * X \ar_{v_2}[r]             & \bchn * X,
  }
\end{align*}
the maps being defined as follows:
\begin{itemize}
  \item For $k \in \N$ let $s_k \colon \map[norm]{\ssim k} X \longrightarrow
        \map[norm]{\ssim k}{\ucov X}$ be a Borel section of the map induced
        by the universal covering map. The existence of such a section 
        is guaranteed by the following theorem -- whose (elementary,
        but rather technical) proof is 
        exiled to the Appendix: 

        \begin{satz}\label{borelsecthm}
          Let $X$ be a connected locally finite \cw-complex or a
          manifold. Then the map 
          \begin{align*}
             P \colon \map[norm]{\ssim k}{\ucov X} & \longrightarrow 
                      \map[norm]{\ssim k}X \\
                      \sigma                 & \longmapsto
                       p \circ \sigma
          \end{align*}
          induced by the universal covering map~$p \colon \ucov X
          \longrightarrow X$ admits a Borel section.
        \end{satz} 

        If $f\in \bfch k X$, we write 
        \begin{align*}
          v_1(f) \colon \map[big]{\ssim k} X & \longrightarrow \R \\
                              \sigma   & \longmapsto
                                         f\bigl( (s_k(\sigma))(e_0)
                                               , \dots, 
                                               (s_k(\sigma))(e_k) 
                                          \bigr).
        \end{align*}
        Since $s_k$ is Borel, $v_1(f)$ is also Borel.   
  \item The map~$v_2$ is given by linear extension
        (cf. Lemma~\ref{mhkrp}).
  \item The map~$v$ is defined in Theorem~\ref{specialthm}: 
        For $f \in \bfch k X$, the homomorphism~$v(f)$ is the linear 
        extension of
        \begin{align*}
          \map[big]{\ssim k} X & \longrightarrow \R \\
                   \sigma & \longmapsto f\bigl( (s_k(\sigma))(e_0)
                                              , \dots 
                                              , (s_k(\sigma))(e_k)
                                         \bigr).
        \end{align*}
\end{itemize}

\begin{qedbem}\label{abde}
  Since the functions living in~$\bfch * X$ are both 
  $\pi_1(X)$-invariant and bounded,
  $v_1$ is a well-defined cochain map, and $H^*(v_1)$ does not
  increase the norm. By construction, the diagram is commutative. In
  particular, $H^*(v_2)$ is surjective by
  Theorem~\ref{specialthm}. 
\end{qedbem}

The crux of the above diagram is that the vertical arrow induces an
isometric isomorphism on the level of cohomology and that $v_1$
does not increase the seminorm. Hence the duality principles allow us
to deduce that the algebraic isomorphism must be isometric for all
connected locally finite \cw-complexes:

\begin{bew}[Proof (of Theorem~\ref{isomisomthm} for connected locally
finite \cw-complexes)] \hfil 
  According to Theorem~\ref{algebraicisomthm}, the induced homomorphism
  $H_*(i_X) \colon \hsing * X \longrightarrow \hmh * X$ is an
  isomorphism. Therefore, it remains to show that~$H_*(i_X)$ is
  compatible with the seminorms.

  Let $k \in \N$ and $\alpha \in \hsing k X$. Since $i_X \colon
  \csing * X \longrightarrow \mhch * X$ is norm preserving, it is
  immediate that
  $\mhnorm{H_k(i_X)(\alpha)} \leq \lone{\alpha}.$

  The proof of the reverse inequality is split into two cases:
\begin{enumerate}
  \item  
  Suppose $\mhnorm[big]{H_k(i_X)(\alpha)} = 0$.  From 
  Lemma~\ref{mhkrp} and Lemma~\ref{mhduality} we obtain
  \[ \krp[big]{H^k(v_2)(\varphi)}{\alpha} 
   = \krp[big]{\varphi}{H_k(i_X)(\alpha)} 
   = 0\]
  for all $\varphi \in \bmc k X$.  By the previous remark, $H^k (v_2)$ 
  is surjective. Hence
  \[ \fa{\psi \in \bch k X} \krp \psi \alpha = 0,\]
  implying $\lone \alpha = 0$ by duality (Theorem~\ref{duality}).
  \item
  Let $\mhnorm[big]{H_k(i_X)(\alpha)} > 0$. 
  In this case, the duality principle for measure homology
  (Lemma~\ref{mhduality}) and Lemma~\ref{mhkrp} yield
  \begin{align*}
    \mhnorm[big]{H_k(i_X)(\alpha)} 
  & \geq \sup
    \Bigl\{ \frac 1 {\supn \varphi} 
    \Bigm| \text{$\varphi \in \bmc k X$, 
                 $\krp[big] \varphi {H_k(i_X)(\alpha)} =1$}
    \Bigr\} \\
  & = \sup
    \Bigl\{ \frac 1 {\supn \varphi} 
    \Bigm| \text{$\varphi \in \bmc k X$,  
                 $\krp[big]{H^k(v_2) (\varphi)} \alpha =1$}
    \Bigr\}.
  \end{align*}
  We compare the last set with the corresponding set of 
  Theorem~\ref{duality}:
  Let $\psi \in \bch k X$ such that $\krp \psi \alpha = 1$. Since
  $H^k(v)$ is an isometric isomorphism (Theorem~\ref{specialthm}),
  there exists a $\xi \in H^k\bigl( \bfch * X \bigr)$ satisfying
  \[ H^k(v) (\xi) = \psi 
     \quad\text{and}\quad \supn \xi = \supn \psi.\]
  Then $\varphi := H^k(A)(\xi) \in \bmc k X$ possesses the following
  properties: 
  \begin{itemize}
    \item By construction,
          $ H^k(v_2) (\varphi) 
           = \bigl(H^k(v_2) \circ H^k(v_1)\bigr)(\xi)
           = H^k(v)(\xi) 
           = \psi,$ 
          and hence 
          \[ \krp[big]{H^k(v_2) (\varphi)}{\alpha} 
           = \krp{\psi}{\alpha} 
           = 1.\]
    \item Furthermore, we get from Remark~\ref{abde} that 
          \[ \supn \varphi  = \supn[big]{H^k(v_1)(\xi)}
             \leq \supn \xi = \supn \psi.\]
  \end{itemize}      
  Combining these properties with the above estimate results in
  \begin{align*} 
    \mhnorm[big]{H_k(i_X)(\alpha)} 
  & \geq \sup
    \Bigl\{ \frac 1 {\supn \psi}
    \Bigm| \text{$\psi \in \bch k X$, $\krp{\psi}{\alpha} = 1$}
    \Bigr\}.
  \end{align*}
  Since $\lone \alpha \geq \mhnorm{H_k(i_X)(\alpha)} > 0$, we can
  use the duality principle (Theorem~\ref{duality}) to conclude that 
  \[\mhnorm{H_k(i_X)(\alpha)} \geq \lone \alpha.\qedhere\]
\end{enumerate}
\end{bew}

\subsection{Second step -- the general case}

We can now reduce the general case of conncected \cw-complexes to the
case of connected finite \cw-complexes:

\begin{bew}[Proof (of Theorem~\ref{isomisomthm} -- the general case)]
  Let $X$ be a connected \cw-complex. 
  Again, Theorem~\ref{algebraicisomthm} states that $H_*(i_X) \colon \hsing
  * X \longrightarrow \hmh * X$ is an isomorphism and it remains to
  prove that $H_*(i_X)$ is isometric:

  Let $\alpha \in \hsing k X$. Then clearly $\lone{\alpha} \geq
  \mhnorm{H_*(i_X)(\alpha)}$. For the converse inequality, 
  let $\mu \in \mhch k X$ be a measure chain
  representing~$H_*(i_X)(\alpha)$. Hansen showed in his proof that
  measure homology respects certain colimits \cite[proof of
  Proposition~5.1]{hansen} that we can find a compact subspace~$A
  \subset X$ and a measure chain~$\nu \in \mhch k A$ such that $\mhch k
  i (\nu) = \mu$, where $i \colon A \hookrightarrow X$ is the
  inclusion. Then, as one can check easily, $\mhnorm \nu = \mhnorm
  \mu$. 

  Since $A$ is compact and $X$ is a connected \cw-complex, we can
  assume that $A$ is a connected finite subcomplex of~$X$. 

  The first step of our proof shows that the isomorphism~$H_*(i_A)
  \colon \hsing * A \longrightarrow \hmh * A$ is isometric. In
  particular, the preimage $\beta := H_*(i_A)^{-1}([\nu])$ satisfies 
  \begin{align*}
    \lone\beta & = \mhnorm{[\nu]} \leq \mhnorm{\nu} = \mhnorm{\mu}. 
  \end{align*}
  By construction, 
  $\hsing * i (\beta) 
    = H_*(i_X)^{-1}\circ \hmh * i \circ H_*(i_A)(\beta)
    = \alpha$, and 
  therefore functoriality implies
  \[ \lone{\alpha} =    \lone[big]{H_*(i)(\beta)}
                   \leq \lone{\beta}
                   \leq \mhnorm{\mu}.\]
  Taking the infimum over all representatives~$\mu$
  of~$H_*(i_X)(\alpha)$ gives the desired inequality~$\lone{\alpha}
  \leq \mhnorm{H_*(i_X)(\alpha)}$.
\end{bew}

\section{A glimpse at the smooth universe}\label{smoothsec}

In this section, a short exposition of the smooth version of the
isometric isomorphism (Theorem~\ref{smisomisomthm}) is given. We first
state a precise definition of smooth measure
homology. Section~\ref{smsingsec} introduces smooth singular homology
which is the building bridge between singular homology and smooth
measure homology. In
Section~\ref{smalgisomsec} and~\ref{smisomisomsec}, the corresponding
algebraic and isometric isomorphisms are explained. 

\subsection{Definition of smooth measure homology}

In order to define smooth measure homology, we have to make precise
what smooth simplices are and what the topology on the corresponding
mapping spaces looks like. Then the definition is completely analogous
to the definition of measure homology:

\begin{defi}
  Let $M$ be a smooth manifold and let $k \in \N$. 
  \begin{itemize}
    \item A singular simplex~$\sigma \colon \ssim k \longrightarrow M$ 
          is called \defin{smooth} if it can be extended to a smooth
          map on an open neighbourhood of~$\ssim k$.  
          We write $\maps[norm]{\ssim k}{M}$ for the set of all smooth 
          singular simplices.
    \item The \cone\defin{-topology} on~$\maps[norm]{\ssim k} M$ is
          the unique topology that turns the differential
          $\maps[norm]{\ssim k}M \longrightarrow \map[norm]{T\ssim
          k}{TM}$ into a homeomorphism onto the image, where
          $\map[norm]{T\ssim k}{TM}$ is endowed with the compact-open 
          topology.
    \item The $k$\defin{-th smooth measure chain group}~$\smmhch k M$
          is defined like the $k$-th measure chain group but using
          $\maps[norm]{\ssim k}M$ with the \cone-topology instead of
          $\map[norm]{\ssim k}M$ with the compact-open topology. 
    \item The $k$\defin{-th smooth measure homology group}~$\hsmmh k
          M$ is the $k$-th homology group of~$\smmhch * M$, where the
          boundary operator is defined as in the non-smooth case. 
    \item The total variation on~$\smmhch * M$ induces a seminorm on
          smooth measure homology, which is denoted
          by~$\smhnorm{\,\cdot\,}.$
          \qedhere 
  \end{itemize}
\end{defi}

\subsection{Smooth singular (co)homology}\label{smsingsec}

There is no reasonable chain map between~$\csing * M$ and~$\smmhch *
M$, so we take a small detour to smooth singular homology. On
(co)homology, it turns out -- as one would suspect -- that smooth
singular homology and singular homology are isometrically isomorphic. 

\begin{defi}
  Let $M$ be a smooth manifold. 
  Then $\smsing * M$ stands for the
  subcomplex of~$\csing * M$ generated by all smooth simplices. We
  write 
  \[\hsmsing * M := H_* \bigl( \smsing * M 
                            , \bou|_{\smsing * M}
                       \bigr)\]
  for \defin{smooth singular homology} (with real
  coefficients). Furthermore, we obtain a seminorm
  on~$\hsmsing * M$ induced by the $\ell^1$-norm
  on~$\smsing * M$. 
\end{defi}

\begin{prop}\label{smsingisomprop}
  Let~$M$ be a smooth manifold. The inclusion~$j_M \colon \smsing * M
  \hookrightarrow \csing * M$ induces a natural (isometric)
  isomorphism
  \[ \hsmsing * M \cong \hsing * M.\]
\end{prop}
\begin{bew}
  Via the Whitney approximation theorem, a smoothing operator
  \[ s \colon \csing * M \longrightarrow \smsing * M \]
  can be constructed \cite[page~417]{lee}, satisfying the following
  conditions: The map~$s$ is a chain map with $s \circ j_M = \id$ and
  $j_M \circ s \simeq \id$, and for each singular simplex~$\sigma \in
  \map[norm]{\ssim k}M$ the image $s(\sigma) \in \smsing k M$
  consists of just one smooth simplex.

  The first part implies that $H_*(j_M)  \colon \hsmsing * M
  \longrightarrow \hsing * M$ is an isomorphism with
  inverse~$H_*(s)$. Moreover, the second property ensures that
  $H_*(s)$ does not increase the seminorm. Hence $H_*(j_M)$ is
  isometric. 
\end{bew}

Moreover, we need a smooth version of bounded cohomology:

\begin{defi}
  If $M$ is a smooth manifold, we define~$\bsmcch * M$ as the set of
  all homomorphisms~$f \colon \smsing k M \longrightarrow
  \R$ satisfying
  \[ \sup\bigl\{ |f(\sigma)|
         \bigm| \sigma \in \maps[norm]{\ssim k}{M}
         \bigr\} 
     < \infty. \]
  As in the non-smooth case, $\bsmcch * M$ can be equipped with a
  coboundary operator~$\cobou$ and we write   
  $\bsmc * M := H^*\bigl( \bsmcch * M
                         , \cobou
                    \bigr)$
  for \defin{smooth bounded cohomology}.
\end{defi}

\begin{prop}\label{bsmsingisomprop}
  Let $M$ be a smooth manifold. Then the restriction homomorphism
  $\bchn * M \longrightarrow \bsmcch * M$ induces a natural isometric
  isomorphism 
  \[ \bsmc * M \cong \bch * M.\]
  Moreover, there is a duality principle (in the sense of
  Theorem~\ref{duality}) for smooth bounded cohomology.
\end{prop}
\begin{bew}
  The dual 
  \begin{align*}
    \bchn * M & \longrightarrow \bchn {*-1} M \\
            f & \longmapsto     (-1)^* \cdot f \circ h_{*-1}
  \end{align*}
  of a (bounded) chain homotopy~$h \colon j_M \circ s \simeq \id$ shows that
  $\bsmc * M \cong \bch * M$. 
  Note that $h$ can be chosen to be bounded in each degree
  \cite[page~417ff]{lee}, so that the above cochain homotopy is indeed 
  well-defined. 
  Furthermore, the duals $\bchn * {j_M}$
  and $\bchn * s$ do not increase the seminorm. Therefore, $\bch *
  {j_M}$ is an isometric isomorphism.

  The duality principle follows easily from Theorem~\ref{duality} and
  the fact that both $H_*(j_M)$ and~$\bch * {j_M}$ are isometric
  isomorphisms. 
\end{bew}

\subsection{The algebraic isomorphism}\label{smalgisomsec}

Analogously to the non-smooth case, we can compare smooth singular
homology and smooth measure homology:

\begin{defi}
  If $M$ is a smooth manifold and $k \in \N$, we write
  \begin{align*}
    \sminc_M \colon \smsing k M & \longrightarrow \smmhch k M \\
    \sum_{\sigma \in \maps[norm]{\ssim k}M} a_\sigma \cdot \sigma
    & \longmapsto
    \sum_{\sigma \in \maps[norm]{\ssim k}M} a_\sigma \cdot \delta_\sigma.
    \qedhere
  \end{align*}
\end{defi}

Clearly, this inclusion induces a natural cochain map $\smsing * M
\longrightarrow \smmhch * M$ (also denoted by~$\sminc_M$) which is
norm-preserving. 

\begin{satz}
  Let $M$ be a smooth manifold. Then $\sminc_M$ 
  induces a natural isomorphism (of real vector spaces)
  \[ \hsmsing * M \cong \hsmmh * M.\] 
\end{satz}
\begin{bew}
  Zastrow explains how one can translate his proofs of the
  Eilenberg-Steen\-rod axioms for measure homology to the smooth case
  \cite[Theorem~3.4]{zastrow}. Then an ``induction'' similar to Milnor's
  proof of Poincar\'e duality \cite[page~351,
  Theorem~(4.10)]{massey,strohm} shows that $\sminc_M$ induces an
  isomorphism between smooth singular homology and smooth measure homology. 
\end{bew}

\begin{koro}
  For all smooth manifolds, singular homology and smooth measure
  homology are naturally isomorphic. 
  \hfill $\square$
\end{koro}

\subsection{The isometric isomorphism}\label{smisomisomsec}

As in Section~\ref{firststep}, we can apply duality to see that smooth
singular homology and smooth measure homology are isometrically
isomorphic. 

The duals $\bsmmhch * M$ and~$\bcsmmh * M$ are defined like $\bmcch *
M$ and~$\bmc * M$ where the mapping space 
$\map[norm]{\ssim k}M$ (equipped with the
compact-open topology) is replaced by~$\maps[norm]{\ssim k}M$ (with
the \cone-topology). Literally the same proof as for Theorem~\ref{mhduality}
yields that there is a duality principle for $\bcsmmh * M$ and~$\hsmmh
* M$. 

We now have collected all the tools necessary to prove
Theorem~\ref{smisomisomthm}: 

\begin{satz}
  Let $M$ be a connected smooth manifold. Then the natural isomorphism
  $H_*(\sminc_M) \circ H_*(j_M)^{-1}
     \colon 
     \hsing * M    
     \longrightarrow
     \hsmmh * M
  $
  induced by the canonical inclusions $\sminc_M \colon \smsing * M
  \longrightarrow \smmhch * M$ and~$j_M \colon \smsing * M
  \longrightarrow \csing * M$ is isometric.  
\end{satz}
\begin{bew}
  Analogously to Section~\ref{firststep}, we can consider the commutative
  diamond
  \begin{equation*}
   \xymatrix@=3em@dr{%
    \bfch  * M  \ar^{v}[r] \ar_{v_1}[d] 
  & \bchn * M,   \ar^{v_3}[d]\\
    \bsmmhch * M  \ar_{v_2}[r]        
  & \bsmcch * M 
    }
  \end{equation*}
  where $v_1$, $v_2$ and $v$ are defined as in
  Section~\ref{firststep}. The cochain map~$v_3$ is given
  by restriction, which is an isometric isomorphism by
  Proposition~\ref{bsmsingisomprop}. 
  Since we have established duality principles for smooth singular
  homology and  smooth measure homology, 
  the same proof as in Section~\ref{firststep} shows
  that $H_*(\sminc_M)$ is an isometric isomorphism. Now the theorem
  follows with help of Proposition~\ref{smsingisomprop}.
\end{bew}

\section{Applications}\label{applications}

We can apply the isometric isomorphisms of Theorem~\ref{isomisomthm}
and~\ref{smisomisomthm} to compute the simplicial volume in terms of
(smooth) measure homology. In fact, this was the motivation for
Thurston to invent measure homology. 

\begin{defi}
  If $M$ is an oriented closed connected manifold, its
  \defin{simplicial volume} is given by
  \[ \sv{M} := \lone{\fclr M},\]
  where $\fclr M \in \hsing{\dim M} M$ denotes the (real) fundamental
  class of~$M$.  
\end{defi}

The simplicial volume of spheres and tori is zero since these
manifolds admit selfmaps of degree larger than~1. On the other hand,
the simplicial volume of hyperbolic manifolds is non-zero
\cite[Theorem~C.4.2]{bp}, e.g., $\sv{F_g} = 4g-4$ for all oriented
closed surfaces~$F_g$ of genus~$g\geq 2$. An amazing facet of
the simplicial volume is that it is a \emph{topological} invariant
bounding the minimal volume from below (modulo a constant factor
depending only on the dimension) \cite[page~12]{gromov}. 

\begin{koro}
  Let $M$ be an oriented closed connected manifold. Then
  \[ \sv{M} = \mhnorm{H_*(i_M)(\fclr M)}.\]
  If $M$ is smooth, then
  \[ \sv{M} = \smhnorm{H_*(\sminc_M) \circ H_*(j_M)^{-1}      
                       (\fclr M)}.\]
\end{koro}
\begin{bew}
  The first part follows because all closed manifolds are \cw-complexes
  (up to homotopy) 
  and hence Theorem~\ref{isomisomthm} is applicable. The second part
  is a direct consequence of Theorem~\ref{smisomisomthm}.
\end{bew}

For hyperbolic manifolds, this result is well-known
\cite[Remark~0.1]{zastrow}. 

The above corollary makes it possible to apply Thurston's smearing
technique \cite[page~6.8ff, page~547ff,
Chapter~5]{thurston,ratcliffe,strohm} to prove the proportionality principle of
simplicial volume \cite[page~6.9]{thurston} (more details are given in
\cite[Chapter~5]{strohm}):

\begin{satz}[Proportionality Principle of Simplicial Volume]
  Let $M$ and~$N$ be Riemannian manifolds with isometric universal
  Riemannian coverings. Then
  \[ \frac{\sv M}{\vol M} = \frac{\sv N}{\vol N}.\]
\end{satz}

Gromov's original proof of the proportionality principle is sketched
in \emph{Volume and bounded cohomology}~\cite[page~11]{gromov}.

\begin{koro}
  The simplicial volume of oriented closed connected smooth flat
  manifolds vanishes. 
\end{koro}
\begin{bew}
  If we scale a flat manifold, the simplicial volume and the
  universal Riemannian covering space remain the same but the volume
  changes. Hence the proportionality principle implies that the
  simplicial volume must be zero. 

  Note that this result can also be obtained by Gromov's estimate of
  the minimal volume \cite[page~12]{gromov} or by the boundedness
  results for the Euler class of Milnor and Sullivan \cite[page~23]{gromov}.
\end{bew}

Another consequence of the proportionality principle is the following
mapping theorem, due to Gromov \cite[5.36]{gromovb}. 

\begin{koro}[Gromov's mapping theorem] Let $n \in \N$ and let $S_1,
  \dots, S_n$ and $S'_1, \dots, S'_n$ be hyperbolic surfaces. If $f
  \colon S_1 \times \dots \times S_n \longrightarrow S'_1 \times \dots \times
  S'_n$ is a continuous map of degree~$d$, then
  \[     |\chi(S_1 \times \dots \times S_n)| 
    \geq d \cdot |\chi(S'_1 \times \dots \times S'_n)|.\]
\end{koro}
\begin{bew}
  For brevity, we write $S:= S_1 \times \dots \times S_n$ and $S' :=
  S'_1 \times \dots \times S'_n$. 
  Since the surfaces involved are all hyperbolic, the universal
  Riemannian coverings of~$S$ and~$S'$ coincide. Therefore, we obtain
  from the proportionality principle and by functoriality
  of~$\lone{\,\cdot\,}$ that 
  \[      \frac{\sv{S'}}{\vol S'}
     =    \frac{\sv S}{\vol S}
     \geq \frac{d \cdot \sv{S'}}{\vol S}. \]
  Moreover, $\sv{S'} \neq 0$ \cite[example on page~9 and (1) on
  page~10]{gromov}, and both the volume
  and the Euler characteristic are multiplicative with respect to the
  product. Since the hyperbolic volume and the absolute value of the
  Euler characteristic of oriented closed connected surfaces are
  proportional (by Gau\ss-Bonnet), the result follows.
\end{bew}

\setcounter{section}{0}
\renewcommand{\thesatz}{(A.\arabic{satz})}
\setcounter{satz}{0}
\section*{Appendix -- Existence of a Borel section}\label{appendix}

To complete the proof of Theorem~\ref{isomisomthm}, we still have to
provide a proof of Theorem~\ref{borelsecthm}. As a first step we 
prove the following (stronger) statement:

\begin{lem}\label{borelseclem}
  Let $X$ be a locally path-connected, semi-locally simply
  connected space such that the universal covering~$\ucov X$ is
  metrisable  
  (e.g., let $X$ be a locally finite \cw-complex). Then the map  
  \begin{align*}
     P \colon \map[norm]{\ssim k}{\ucov X} & \longrightarrow 
              \map[norm]{\ssim k}X \\
              \sigma                 & \longmapsto
              p \circ \sigma
  \end{align*}
  induced by the universal covering map~$p \colon \ucov X
  \longrightarrow X$ is a local homeomorphism. 
\end{lem}

In the proof of this lemma we use the following notation
for the sub-basic sets of the compact-open topology:

\begin{defi}
  In the above situation, if $K \subset \ssim k$ is compact and $U
  \subset X$ (or $U \subset \ucov X$) is open, we write $U^K$ for the
  set of all~$f \in \map[norm]{\ssim k}{X}$ (or all~$f \in
  \map[norm]{\ssim k}{\ucov X}$ respectively) satisfying~$f(K) \subset U$.
\end{defi}

\begin{bew}
Let $\sigma \in \map[norm]{\ssim k}{\ucov X}$. Then there is a small
neighbourhood~$U$ of~$\sigma(e_0)$ in~$\ucov X$:  

  \begin{defi}
    An open subset~$U \subset \ucov X$ is
    called \defin{small} if $p$ is
    trivial over~$p(U)$ and the restriction
    $p|_{U} \colon U \longrightarrow p(U)$
    is a homeomorphism.  
  \end{defi}

In particular, the image~$p(U) \subset X$ is open because covering maps are
open. 

  \begin{bem}\label{smallrem}
    Since $p \colon \ucov X \longrightarrow X$ is a covering map and
    since $X$ is locally path-connected, each point in~$\ucov X$ 
    possesses a basic family of small neighbourhoods. 
  \end{bem}

We show that $P( U^{\{e_0\}})$ is open and that the restriction 
$P|_{U^{\{e_0\}}} \colon U^{\{e_0\}} \longrightarrow P(U^{\{e_0\}})$
is a homeomorphism.

\emph{The set $P(U^{\{ e_0\}})$ is open in~$\map[norm]{\ssim k} X$.} By
  definition, $P(U^{\{ e_0\}}) \subset \bigl(p(U)\bigr)^{\{e_0\}}$. On
  the other hand, for each $\tau \in \bigl(p(U)\bigr)^{\{e_0\}}$ there
  exists a lift $\ucov \tau \colon \ssim k \longrightarrow \ucov X$
  such that $\ucov \tau(e_0) \in U$ because $\ssim k$ is simply
  connected. 
  Thus
  \[ P(U^{\{ e_0\}}) = \bigl(p(U)\bigr)^{\{e_0\}}.\]
  Since $p(U)$ is open in $X$, this is an open subset of~$\map[norm]{\ssim k}
  X$. 

  \emph{The restriction $P|_{U^{\{e_0\}}} \colon U^{\{e_0\}}
  \longrightarrow \bigl( p(U)\bigr)^{\{e_0\}} =
  P\bigl(U^{\{e_0\}}\bigr)$ is bijective.}
  Since $U$ is small, $p|_U$ is injective. Hence the uniqueness of
  lifts (prescribed on $e_0$ by the property to map into $U$) shows 
  injectivity of~$P|_{U^{\{e_0\}}}$ ($\ssim k$ is connected). 

  \emph{The restriction $P|_{U^{\{e_0\}}} \colon U^{\{e_0\}}
  \longrightarrow \bigl( p(U)\bigr)^{\{e_0\}} =
  P\bigl(U^{\{e_0\}}\bigr)$ is a homeomorphism.} 
  By definition of the compact-open topology, $P$ is continuous. It therefore
  remains to prove that the restriction $P|_{U^{\{e_0\}}}$ is open
  (which is the lion share of the proof): 

  Let $\ucov\tau \in U^{\{ e_0\}}$ and let $A$ be an open neighbourhood
  of~$\ucov\tau$ in $U^{\{e_0\}}$. 
  We have to show that $P|_{U^{\{e_0\}}}(A) \subset \map[norm]{\ssim
  k}{X}$ is open: 
  Since $P|_{U^{\{e_0\}}}$ is injective, this restriction is
  compatible with unions and intersections. By definition of the
  compact-open topology, it is therefore sufficient to consider the
  case $A = V^K \cap U^{\{e_0\}}$, where $K \subset \ssim k$ is
  compact and $V \subset \ucov X$ is open.

  In the following, we make use of especially small subsets
  of~$\ucov X$:

  \begin{defi}
    A family~$(U_i)_{i \in I}$ of subsets of $\ucov X$ is
    \defin{tiny} if the following conditions hold: all $U_i$ are small
    and whenever $U_i \cap U_j \neq \emptyset$, then the union~$U_i \cup
    U_j$ is also small.
  \end{defi}

  \begin{prop}
    Let $Y \subset \ucov X$ a compact subset and
    let $(U_i)_{i \in I}$ be a family of open subsets of~$\ucov X$
    covering~$Y$. Then there is a number~$\varepsilon \in \R_{>0}$
    with the following property:
    If $x \in Y$ and $B_\varepsilon(x)$ is the open ball in~$\ucov X$
    around~$x$ with radius~$\varepsilon$, then $B_\varepsilon(x)
    \subset U_i$ for some $i \in I$. 
    Such an~$\varepsilon$ is called a \defin{Lebesgue number} of the 
    family~$(U_i)_{i \in I}$.
  \end{prop}
  \begin{bew}
    One can use literally the same proof as in Dugundji's book for
    the existence of a Lebesgue number (in a slightly weaker
    context) \cite[Theorem~XI.4.5]{dugundji}.
  \end{bew}

  \begin{koro}\label{tinykoro}
    Let $Y \subset \ucov X$ be a compact subset covered by a
    family~$(V_j)_{j \in J}$ of open sets. Then there is a (finite)
    tiny family~$(U_i)_{i \in I}$ covering~$Y$ subordinate
    to~$(V_j)_{j \in J}$.
  \end{koro}
  \begin{bew}
    It is possible to cover~$Y$ by a
    family~$(V'_j)_{j \in J'}$ of small sets subordinate to~$(V_j)_{j
    \in J}$ (Remark~\ref{smallrem}). 
    Let $\varepsilon \in \R_{>0}$ be a Lebesgue-number, in the above
    sense, of this covering. 
    Since $\ucov X$ is locally path-connected,  
    there is a covering~$(U_i)_{i \in I}$
    of~$Y$ by small subsets satisfying
    \[ \diam(U_i) < \frac{\varepsilon}2\]
    for all~$i \in I$. Then $(U_i)_{i \in I}$ is tiny: Let $i$,~$j \in
    I$ with $U_i \cap U_j \neq \emptyset$. Thus
    \[ \diam(U_i \cup U_j) < \varepsilon. \]
    By construction of~$\varepsilon$, there is an $\ell \in J$ such that
    $U_i \cup U_j \subset V'_\ell$. 
    Since $U_i$ and~$U_j$ are open, so is their union~$U_i \cup
    U_j$. Now $U_i \cup
    U_j \subset V'_\ell$ implies that $U_i \cup U_j$ is small. 
    Furthermore, we can choose $I$ to be finite because $Y$ is
    compact. 
  \end{bew}

  Using the above corollary, we obtain a tiny covering~$(V_j)_{j \in
  J}$ of~$\ucov \tau(\ssim k)$. Applying the above corollary 
  twice more (on the compact sets~$\ucov \tau(K)$ and $\ucov
  \tau(e_0)$ and the induced coverings~$(V \cap V_j)_{j \in J}$ and
  $(U \cap V_j)_{j \in J}$), we can find a
  finite tiny covering~$(U_i)_{i \in I}$ of~$\ucov \tau (\ssim k)$ and
  compact subsets~$(K_i)_{i \in I}$ of~$\ssim k$ such that the
  intersection
  \[ \ucov B := \bigcap_{i \in I}U_i^{K_i}\]
  satisfies
  \[ \ucov B \subset U^{\{ e_0 \}} \quad\text{and}\quad
     \ucov B \subset V^K. \]

  By construction, $\ucov B$ is open in~$\map[norm]{\ssim k}{\ucov X}$
  and $\ucov B \subset V^K \cap U^{\{e_0\}}= A$. 
  It therefore suffices to show that
  $P|_{U^{\{ e_0\}}}(\ucov B) \subset \map[norm]{\ssim k}X$ is
  open. More precisely, we prove that
  \[ P|_{U^{\{e_0\}}}(\ucov B) = B,\]
  where
  \[ B := \bigcap_{i \in I} \bigl( p(U_i)
                            \bigr)^{K_i}. \]
  It is clear that $P|_{U^{\{e_0\}}}(\ucov B) \subset B$. Conversely,
  let $\rho \in B$. It suffices to check that the unique lift~$\ucov
  \rho \in \map[norm]{\ssim k}{\ucov X}$ with $\ucov \rho(e_0) \in U$
  lies in~$\ucov B$. In the following, we prove that the set
  \[ D := \{ x \in \ssim k
          \mid \fa{i \in I(x)} \ucov \rho(x) \in U_i
          \}\]
  is open and closed and that it contains~$e_0$, where we used the notation
  \begin{align*}
    I(x)       & := \{ i \in I \mid x \in K_i \}. 
  \end{align*}

  The key to proving this claim is the following lemma based on tininess: 

  \begin{lem}\label{ilemma}
    If $x \in \ssim k$ and if there is a~$j \in I(x)$ such that
    $\ucov\rho(x) \in U_j$, then $\ucov \rho(x) \in U_i$ for
    \emph{all}~$i \in I(x)$.
  \end{lem}
  \begin{bew}
    Let $i \in I(x)$. Because of $\ucov \tau(K_i) \subset U_i$ and $\ucov
    \tau(K_j) \subset U_j$, we obtain~$U_i \cap U_j \neq \emptyset$.
    Hence $U_i \cup U_j$ is small (the family~$(U_\ell)_{\ell \in I}$
    is tiny). In particular, 
    $p^{-1} \bigl(\rho(x)\bigr) \cap (U_i \cup U_j)$
    contains precisely one element (namely $\ucov \rho(x) \in
    U_j$). But $i \in I(x)$ implies that
    $p^{-1}\bigl(\rho(x)\bigr) \cap U_i$
    also has to contain (exactly) one element. Therefore, $\ucov \rho
    (x) \in U_i \cap U_j$, which shows $\ucov \rho(x) \in U_i$, as
    desired. 
  \end{bew}

  \begin{itemize}
    \item \emph{The set $D$ is open:} For each~$x \in D$,
          \[ W :=   \bigcap_{i \in I(x)} \ucov\rho^{-1}(U_i)
               \cap \bigcap_{i \in I\setminus I(x)} (\ssim k \setminus K_i) \]
          is an open subset of~$\ssim k$ with $x \in W$ and $W \subset D$.
    \item \emph{The set~$D$ is closed:} Let $x \in D\setminus \ssim
          k$. By the above lemma, $\ucov \rho(x) \not\in U_i$ for
          all~$i \in I(x)$. Since $\ssim k = \bigcup_{i \in I}
          \overset \circ K_i$, there is an~$i \in I$ such that $x \in
          \overset\circ K_i$ and $\ucov \rho \not \in U_i$. Let
          \[ W := \overset\circ K_i 
                  \cap 
                  \ucov\rho^{-1} \bigl(p^{-1}(p(U_i)) \setminus U_i
                                 \bigr). \]
          Since $p$ is trivial over~$p(U_i)$ (with discrete fibre),
          the preimage
          \[\ucov\rho^{-1}\bigl(p^{-1}(p(U_i)) \setminus U_i\bigr)\] 
          is open. By construction, $x \in W$ and $W
          \subset \ssim k \setminus D$. Therefore, $D$ is closed.
    \item \emph{The vertex~$e_0$ lies in~$D$:} This is a direct
          consequence of~$\rho \in B \subset p(U)^{\{e_0\}}$ and
          Lemma~\ref{ilemma}.  
  \end{itemize}
  Since $\ssim k$ is connected, this implies $D=\ssim k$ and hence 
  proves Lemma~\ref{borelseclem}.
\end{bew}

Finally, we are able to conclude that $P$ possesses a Borel section,
as claimed in Theorem~\ref{borelsecthm}:

\begin{bew}[Proof (of Theorem~\ref{borelsecthm})]
  We can cover $\map[norm]{\ssim k}{\ucov X}$ with countably many open
  sets~$(V_n)_{n\in\N}$ on which $P$ is a homeomorphism and such that
  $P(V_n)$ is open since $\ucov X$ is second countable (e.g., one
  could take a countable covering of~$\ucov X$ by small sets~$U$ and
  consider sets of the form $U^{\{e_0\}}$). Setting
  $W_0 := P(V_0)$ and
  \[ \fa{n \in\N} W_{n+1} := P(V_{n+1}) \setminus 
                             \bigcup_{%
                             {{j \in \{0, \dots, n\}}}} W_j,
  \]
  we get a countable family $(W_n)_{n\in\N}$ of mutually disjoint
  Borel sets in 
  $\map[norm]{\ssim k}{X}$ such that the inverse $P^{-1}|_{W_n}$ is 
  well-defined
  and continuous for each~$n \in \N$. 
  Moreover, $\map[norm]{\ssim k}X$ is covered by the $(W_n)_{n \in \N}$
  because $P$ is surjective. 
  Putting all these maps together
  yields the desired Borel section of~$P$.
\end{bew}


\end{document}